\documentclass[10pt,twoside]{article}

\newtheorem{proposition}{Proposition}[section]

\newtheorem{theorem}[proposition]{Theorem}

\input epsf

\renewcommand{\bar}[1]{\overline{#1}}

\newcommand{\mun}{\mu_n}
\newcommand{\bhc}{\widehat{\beta}_c}

\newcommand{\An}[1]{A_n^{#1}}

\newcommand{\cfcirc}{{{\calF}^{\!\!\!\raise5pt\hbox{\scriptsize o}}}}

\newcommand{\Gan}{\Ga_n}

\newcommand{\Omn}{\Om_n}

\newcommand{\bchat}{\widehat{\beta}_c}

\newcommand{\LDn}{\frac{1}{n^{d-1}}\, \log\, }

\newcommand{\Bn}{{B(n)}}

\newcommand{\cyl}{\, {\rm cyl}\, }

\newcommand{\vol}{\, {\rm vol}\, }

\newcommand{\dist}{\, {\rm dist_{L^1}}}

\newcommand{\DistP}{\, {\rm dist}_P}

\newcommand{\dinf}{\, {\rm d}_{\infty}\, }

\newcommand{\bd}{\partial\, }
\newcommand{\din}{\partial^{\, in}}

\newcommand{\fk}[3]{\Phi^{\, #1}_{#2} \big[ #3 \big] }
\newcommand{\FK}[3]{\Phi^{\, #1}_{#2} \Big[ #3 \Big] }
\newcommand{\fko}[2]{\Phi^{\, #1}_{#2}}

\newcommand{\MU}[3]{\mu ^{\, #1}_{#2} \Big[ #3 \Big] }
\newcommand{\Mu}[3]{\mu ^{\, #1}_{#2} \big[ #3 \big] }
\newcommand{\Muo}[2]{\mu^{\, #1}_{#2}}

\newcommand{\set}[2]{\big\{ #1\, ; \, #2 \big\}}
\newcommand{\Set}[2]{\Big\{ #1\, ; \, #2 \Big\}}

\def\le{\leq}
\def\ge{\geq}

\newcommand{\mb}[1]{\mbox{#1}}
\newcommand{\text}[1]{\mbox{#1}}

\newcommand{\symdif}[2]{#1\, \triangle\, #2}

\newcommand{\rec}[1]{\frac{1}{#1} }

\newcommand{\conn}{\leftrightarrow }

\newcommand{\nempty}{\neq\emptyset}
\renewcommand{\empty}{=\emptyset}

\renewcommand{\limsup}{\mathop{\overline{\rm{lim}}}}
\renewcommand{\liminf}{\mathop{\underline{\rm{lim}}}}

\newcommand{\lsn}{\limsup_{n\rightarrow\infty}}
\newcommand{\lin}{\liminf_{n\rightarrow\infty}}

\newcommand{\limn}{\lim_{n\rightarrow\infty}}

\newcommand{\comment}[1]{}
\newcommand{\kor}[1]{\smash{{#1}^{\!\!\!\raise5pt\hbox{\scriptsize o}}}}

\newcommand{\lb}[1]{\label{#1}}
\newcommand{\bref}[1]{(\ref{#1})}

\newcommand{\be}{\begin{equation}}
\newcommand{\ee}{\end{equation}}
\newcommand{\bea}{\begin{eqnarray}}
\newcommand{\eea}{\end{eqnarray}}
\newcommand{\bean}{\begin{eqnarray*}}
\newcommand{\eean}{\end{eqnarray*}}
\newcommand{\beal}[1]{\begin{eqnarray} & & \lb{#1}\nonumber\\}
\newcommand{\eeal}{\end{eqnarray}}

\newcommand{\de}{\delta}

\newcommand{\eps}{\varepsilon}

\newcommand{\La}{\Lambda}

\newcommand{\Th}{\Theta}

\newcommand{\ga}{\gamma}
\newcommand{\Ga}{\Gamma}

\newcommand{\Om}{\Omega}
\newcommand{\si}{\sigma}

\newcommand{\N}{{\Bbb{N}}}
\newcommand{\R}{{\Bbb{R}}}

\newcommand{\Zd}{{\Bbb{Z}}^d}
\newcommand{\Zdn}{{\Bbb{Z}}^d_n}

\newcommand{\Rd}{{\Bbb{R}}^d}

\newcommand{\calB}{{\cal B}}

\newcommand{\calE}{{\cal E}}
\newcommand{\calF}{{\cal F}}
\newcommand{\calG}{{\cal G}}
\newcommand{\calH}{{\cal H}}
\newcommand{\calI}{{\cal I}}

\newcommand{\calP}{{\cal P}}

\newcommand{\calS}{{\cal S}}

\newcommand{\calU}{{\cal U}}

\newcommand{\calW}{{\cal W}}

\newcommand{\Ecirc}{\smash{{E}^{\!\!\!\!\raise5pt\hbox{\scriptsize o}}}}
\newcommand{\Acirc}{\smash{{A}^{\!\!\raise5pt\hbox{\scriptsize o}}}}
\newcommand{\bro}{\smash{{B}^{\!\!\!\!\raise5pt\hbox{\scriptsize o}}}}
\newcommand{\dro}{\smash{{D}^{\!\!\!\!\raise5pt\hbox{\scriptsize o}}}}
\newcommand{\qro}{\smash{{Q}^{\!\!\!\!\raise5pt\hbox{\scriptsize o}}}}
\newcommand{\lro}{\smash{{L}^{\!\!\!\!\raise5pt\hbox{\scriptsize o}}}}
\newcommand{\uro}{\smash{{U}^{\!\!\!\!\raise5pt\hbox{\scriptsize o}}}}
\newcommand{\ero}{\smash{{E}^{\!\!\!\!\raise5pt\hbox{\scriptsize o}}}}
\newcommand{\cero}{\smash{{\calE}^{\!\!\!\!\raise5pt\hbox{\scriptsize o}}}}
\newcommand{\cgcirc}{\smash{{\calG}^{\!\!\!\raise5pt\hbox{\scriptsize o}}}}
\newcommand{\cdero}{\smash{{\Bbb{E}}^{\!\!\!\!\raise5pt\hbox{\scriptsize o}}}}
\newcommand{\bgcirc}{\smash{{\Bbb{G}}^{\!\!\!\!\raise5pt\hbox{\scriptsize o}}}}
\newcommand{\begcirc}{{{{\Bbb E}\cap\GGG}^{
\!\!\!\!\!\!\!\!\!\raise3pt\hbox{\scriptsize o}}}}

\newcommand{\aro}{\smash{{A}^{\!\!\!\raise5pt\hbox{\scriptsize o}}}}
\newcommand{\thro}{\smash{{\Th}^{\!\!\!\!\!\raise5.2pt\hbox{\scriptsize o}}}}

\newcommand{\GGG}{{\Bbb{G}}}

\newcommand{\rb}{\partial^{ *}}

\newcommand{\h}{\calH}
\newcommand{\hdu}{\calH^{d-1}}

\input pstricks
\input pst-plot

\def\pointless{\expandafter\PoinTless\the}
{\catcode`p=12 \catcode`t=12 \gdef\PoinTless#1pt{#1}}

\usepackage{graphicx}
\usepackage{amsmath,amsfonts}
\usepackage{Latex-document}

\renewcommand{\thesection}{\arabic{section}}
\renewcommand{\thesubsection}{\arabic{section}.\arabic{subsection}}
\newcommand{\Section}[1]{\section{\hskip -1em.\hskip 1em#1}}
\newcommand{\Subsection}[1]{\subsection{\hskip -1em.\hskip 1em#1}}

\title{\bf Renormalization, Large Deviations \vskip -2mm and Phase Separation in Ising \vskip -2mm and Percolation Models\thanks{The author
thanks the following institutions for their hospitality and/or financial support: Center for Nonlinear Analysis at
CMU, Forschungsinstitut at ETHZ, Universit\'e Paris Sud and IMPA, Brazil. The author has been supported by NSF
Grant DMS-0072217.}\vskip 6mm}

\author{\'Agoston Pisztora\thanks{Department of Mathematical Sciences, Carnegie Mellon University,
Pittsburgh, PA 15213, USA. E-mail: pisztora@andrew.cmu.edu}\vspace*{-0.5cm}}

\date{\vspace{-8mm}}

\begin{document}

\maketitle

\thispagestyle{first} \setcounter{page}{79}

\begin{abstract}

\vskip 3mm

Phase separation is a fairly common physical phenomenon with
examples including the formation of water droplets from humid air
(fog, rain), the separation of a crystalline structure from an
isotropic material such as a liquid or even the formation of the
sizzling gas bubbles when a soda can is opened.

It was recognized long ago (at least on a phenomenological level)
that systems exhibiting several phases in equilibrium can be
described with an appropriate variational principle: the phases
arrange themselves in such a way that the energy associated with
the phase boundaries is minimal. Typically this leads to an almost
deterministic behavior and the phase boundaries are fairly
regular. However, when looked at from a microscopic point of view,
the system consists of a bunch of erratically moving molecules
with relatively strong short-range interaction and the simplicity
of the above macroscopic description looks more than miraculous.
Indeed, when starting from the molecular level, there are many
more questions to be asked and understood: which are the phases
which we will see? why do only those occur? why are the phase
boundaries sharp? how should we find (define) the energy
associated with the interfaces? Only then can we ask the question:
why does the system minimize this energy?

It is only in the last decade that a mathematically satisfactory
understanding of this phenomenon has been achieved. The main goal
of the talk is to present the current state of affairs focusing
thereby on results obtained in joint works with Raphael Cerf. The
connection to fields of mathematics other than probability theory
or statistical mechanics will be highlighted; namely, to geometric
measure theory and to the calculus of variations.

\vskip 4.5mm

\noindent {\bf 2000 Mathematics Subject Classification:} 60K35,
82B, 60F10.

\noindent {\bf Keywords and Phrases:} Ising model, Potts model, Random cluster model, FK percolation, Phase
separation, Phase coexistence, Large deviations, Double bubble, Minimal surfaces.
\end{abstract}

\vskip 12mm

\Section{Introduction}

\vskip-5mm \hspace{5mm}

Although phase separation is a fairly common physical phenomenon
(examples will be given further below) its mathematically
satisfactory understanding, even in the simplest models, has not
been achieved until the last decade. In order to uncover the
mechanism leading to the separation of various phases in an
initially homogeneous material, in particular to explain why and
what kind of phases will occur and which shapes they take, one {
has} to work with a microscopic description of the system, often
at a molecular scale. Materials on such scales however tend to
behave `chaotic' and this strongly motivates (if not compels) the
use of a probabilistic approach. In this approach, the system is
modeled by randomly moving (in equilibrium theory randomly
located) particles which interact with each other according to
some simple, typically short range, mechanism. The goal is then to
derive the large scale behavior of the system based only on the
specification of the local interaction.

The difficulty of the analysis stems from the fact that the
interaction, although local, might be strong enough (depending on
some parameter such as the temperature) to cause a subtile spatial
propagation of stochastic dependence across the entire system. As
a consequence, one has to leave the familiar realm of classical
probability theory whose focus has been laid on the large scale
effects of randomness arising from independent (or weakly
dependent) sources. Instead, we have to deal with a strongly
dependent system and it is exactly this strong dependence which
causes a highly interesting cooperative behavior of the particles
on the macroscopic level which can, in certain cases, be observed
as phase separation.

The problem of phase separation and related issues have been a
driving force behind developing, and a benchmark for testing
various new techniques. Postponing historical remarks until
section 2.1,  let me highlight here only those ones which play an
essential  role in our approach \cite{CP, CP2} achieved in
collaboration with Raphael Cerf. The basic framework is {\em
(abstract) large deviation theory}, see e.g. \cite{Var}, whose
power contributed substantially (admittedly rather to my own
surprize) to the success of this approach. To have sufficient
control of the underlying model, in our case the Ising-Potts
model, we employ {\em spatial renormalization techniques}, as
developed in \cite{Pi}, in conjunction with the Fortuin-Kasteleyn
percolation representation of the Potts model, \cite{FK}. Finally,
tools from geometric measure theory a la Cacciopoli and De Giorgi
will be employed to handle the geometric difficulties associated
with three or higher dimensions, in a similar fashion as was done
in \cite{Ce}.

The goal of this article is to present some recent developments in
the equilibrium theory of coexisting phases in the framework of
the Ising-Potts model. The presentation will be based mainly on
the results contained in \cite{Pi, CP, CP2}. It will include a
description of the underlying physical phenomenon, of the
corresponding mathematical model and its motivation and, following
the statement of the main results, comments on the proofs will be
included. In order to address an audience broader than usual, I
will try to use as little formalism as possible.

\Subsection{Phase separation: examples and phenomenology}

\vskip-5mm \hspace{5mm}

Perhaps the most common and well known example of phase separation
is the development of fog and later rain from humid air. When warm
humid air is cooled down so much that its {\em relative} humidity
at the new cooler temperature would exceed 100\% (i.e. it would
become over-saturated)
 the excess amount of water
precipitates first in the form of very small droplets which we
might observe as fog. The system at this time is not in
equilibrium, rather in a so called metastable state. After waiting
very long time or simply dropping the temperature further down,
the droplets grow bigger and ultimately fall to the ground due to
gravitation in the form of common rain droplets. In this example
phase separation occurred since from a single homogeneous phase
(warm humid air) two new phases have been formed: a cooler mixture
of water and air (note: with 100\% relative
 humidity (saturated) at the new lower temperature) plus a certain amount of
water, more precisely a saturated solution of air in water, in the
form of macroscopic droplets. In fact, in the absence of
gravitation and after very long time, only a huge droplet of fluid
would levitate in a gas (both saturated solutions of air in water
and water in air, respectively). The opposite situation (water
majority, air minority) may also occur. Consider the following
familiar example; think of a bottle of champagne when opened. Here
the change of temperature is replaced by a change of pressure but
the phenomenon is similar with the roles exchanged; at the new
lower pressure the liquid is not able to dissolve the same amonunt
of carbon dioxide, hence this latter precipitates in the form of
small bubbles (droplets), etc.

The phenomenological theory explains this type of phenomenon as
follows. At any temperature there are saturation densities of
air/water and water/air mixtures and only saturated solutions will
coexist in equilibrium. They also determine the volumes of the two
coexisting phases. Moreover the phases arrange themselves in a way
so as to minimize the so-called {\em surface energy}, associated
with the interface between the phases. In fact, it is supposed to
exist a (in general direction-dependent) scalar quantity $\tau$,
called the {\em surface tension}, whose surface integral along the
interface gives the surface energy.  The surface tension, as well
as the saturation densities, have to be measured experimentally.
It is implicitly assumed that the interface is 'surfacelike' and
regular enough so that the integral along the surface makes sense.
The prediction which can be made is that the shape of the phases
in equilibrium is just a solution of the variational principle. By
the classical isoperimetric theorem, in the isotropic case the
solution is just a sphere, hence the occurrence of bubbles and
droplets.  In the non isotropic case the corresponding variational
problem  is called the {\em Wulff problem}. The solution is known
to  be explicitly given \cite{WU, DI, FOb, FOMU, TAa}  by
rescaling appropriately the so called {\em Wulff crystal:}
\[\calW_\tau = \set{x\in\Rd}{ x\cdotp \nu \le \tau(\nu) \mb{ for all
unit vectors }\nu}.
\]
It is worth noting that the same arguments are used to describe
macroscopic crystal shapes as well.

\Subsection{The mathematical model and the goals of the analysis}

\vskip-5mm \hspace{5mm}

The next step is to find a model which is simple enough to be
analyzed by rigorous methods yet rich enough to exhibit the
phenomenon we want to study. In order to accommodate a multitude
of phases we consider a finite number ($q$) possible types of
particles (called colors or spins). Physics suggests to choose a
short range interaction of ``ferromagnetic'' character which means
that particles of identical type prefer to stay together and/or
they repel particles of different types. For simplicity we assume
that the interaction distinguishes only between identical and
different types, otherwise it is invariant under permutation of
colors. There is a standard model of statistical mechanics, called
the (ferromagnetic) {\em  $q$-states Potts model}, which
corresponds exactly to these specifications. We consider the
closed unit cube $\Om\in \Rd$, $d\ge 3$ (modeling the container of
the mixture of particles) overlapped by the rescaled integer
lattice $\Zdn=\Zd/n$. We define $\Omn = \Om \cap \Zdn$, and denote
by $\din \Omn$ the internal vertex boundary of $\Omn$. At each
lattice point $x$ there is a unique particle $\si_x$ of one of the
types $1,2,...,q$. The energy $H(\si)$  of a configuration
$\si=(\si_x)_{x\in\Omn}$ can be chosen to be the number of nearest
neighbor pairs of different types of particles corresponding to
nearest neighbor repulsion. According to the Gibbs formula, the
probability of observing a configuration $\si$ is proportional to
$e^{-\beta H(\si)}$, where $\beta = 1/T$ is the `inverse
temperature' which adjusts the interaction strength. (High $T$ =
large disorder = {\em relatively} small interaction, etc.) Note
that we use a static description of the equilibrium system. It
corresponds to a snap-shot of the system at a given time and the
task is to understand the `typical' picture we will see.

A {\em restricted ensemble} is a collection of certain feasible
configurations. For instance, in the situation of the water/air
mixture every configuration with a fixed number of water and air
particles is possible, and this collection forms our restricted
ensemble. It turns out that the direct study of this particular
ensemble is extremely difficult and it is a crucial idea
(discovered long ago) to go over to a larger, more natural
ensemble, namely to that without any restrictions on the particle
numbers. Then, the restricted system can be regarded as a very
rare event $=: G_n$ in the large ensemble and conditional
probabilities can be used to describe the restricted system. The
events $G_n$ are often in the large deviations regime, and it is
from here that large deviations theory enters the analysis in an
essential way.

The unrestricted system is usually referred to as the Potts model
with free boundary conditions. In the case $q=2$, it is equivalent
to the classical Ising model. The Gibbs formula and the energy
uniquely determines the probability measure in this (and in every)
ensemble. We can introduce mixed boundary conditions as follows.
Divide the boundary $\Ga =
\partial\Om$ of the 'container'
into $q+1$ parts indexed by $\Ga^0,\Ga^1,....\Ga^q$. The parts can
be fairly general but the $(d-1)$-dimensional Hausdorff measure of
their relative boundaries has to be zero. We set for $n\in\N$ and
$i=0,\dots,q$, \bean \Ga^i_n &=& \set{ x\in \Gan}{\dinf(x, \Ga^i)
< 1/n \,\text{ and}\,\ \forall j < i, \ \dinf(x,\Ga^j)\ge 1/n}\ \
\ \ i=0,\dots,q \eean where $\dinf$ denotes the distance
corresponding to the $\max$-norm. We use the sequence of
$q+1$-tuples of sets $\ga(n) = (\Ga^0_n,\dots,\Ga^q_n)$ to specify
boundary conditions by imagining that all particles in $\Ga^j_n$
are of type $j$ for $j=1,...,q$ and none occupies $\Ga^0_n$. This
defines a restricted ensemble and the corresponding
 probability measure is denoted by $\mun = \mun^{\ga(n), \beta,
q}$. The choice  of b.c.s $\ga(n)$ is understood to be fixed.
Let's consider the Ising-Potts model in an $n\times n$ lattice box
$B(n)$ with boundary conditions $j\in\{1,2, ...,q\}$ at a fixed
inverse temperature $\beta$ with the corresponding probability
measure $\Muo{(j),\beta}{\Bn}$. It is well known that, as
$n\to\infty$, a unique, translation invariant infinite volume
measure $\Muo{(j), \beta}{\infty}$ emerges as the weak limit of
the sequence $(\Muo{(j),\beta}{\Bn})_{n\ge 1}$. We can define the
{\em order parameter} $\theta = \theta(\beta, q,d)$ as the {\em
excess density} of the dominant color, the excess is measured from
the symmetric value $1/q$. The model exhibits phase transition in
the sense that in dimensions $d\ge 2$ there exists a critical
value $0< \beta_c(d)  < \infty$ such that for $\beta < \beta_c$,
$\ \theta(\beta) := \Mu{(j), \beta}{\infty}{\si_0 = j} - 1/q = 0$
but for $\beta > \beta_c$, $\theta(\beta) > 0$, i.e., when the
interaction becomes strong enough, the influence of the
(arbitrarily) far away boundary still propagates all the way
through the inside of the volume and creates a majority of
$j$-type particles. (Note that in the Ising model ($q=2$), the
spontaneous magnetization $m^*$ is equal to $\theta$.) The
probability measures $\Muo{(j), \beta}{\infty}$, $j=1,2,...,q$,
describe in mathematical terms what we call 'pure' phases and
which correspond to the saturated solutions in the initial
example.

Having chosen our model, let us formulate the goals of the
analysis. Clearly, the main goal is to verify the predictions of
the phenomenological theory, namely, that on the macroscopic scale
the phases will be arranged according to some solution of the
variational principle corresponding to minimal surface energy.
First, however, the participating phases have to be found and
identified. Moreover, as pointed out earlier, the previous
statement contains a couple of implicit assumptions, such as the
existence of the surface tension, the absence of transitional
states (where one phase would smoothly go over into another one)
the regularity of the interface boundaries, etc., all of which
have to be justified from a microscopic point of view.

\Subsection{Connection to minimal surfaces}

\vskip-5mm \hspace{5mm}

In this section a partially informal discussion of some examples
will be presented with the aim to make the close relation to
minimal surfaces transparent.
\medskip
%\begin{figure}[h]
\vskip-40pt \centerline{ \hskip-30pt \hbox{ \psset{unit=0.1cm}
\pspicture(-25,-12)(40,40) \rput{0}(-5,0){ \rput{0}(0,0){
\newdimen\X
\newdimen\Y
\newdimen\C
% globalement on avance de 18unites
%\X=0.8cm
%\Y=-1.6cm
\X=0.4cm \Y=-.8cm \C=0.8pt
%\newcount\couleur
%\global\couleur=0.1
\psdots(21,7)
%\psline[linewidth=1.5pt]{-}(30,11)(21,7)(12,3)
%\psline[linestyle=dashed,linewidth=1.5pt]{-}(12,3)(21,5)(30,11)
\psline[linestyle=dashed,linewidth=1.5pt]{-}(18,4)(21,7)(18,22)
\psline[linewidth=1.5pt]{-}(36,4)(21,7)(36,22)
\psline[linewidth=1.5pt]{-}(24,-8)(21,7)(24,10)
\psline[linewidth=1.5pt]{-}(6,-8)(21,7)(6,10)
%\psline[linewidth=2pt]{-}(6,1)(24,1)(36,13)(18,13)(6,1)
% full lines upper half
\psline[linewidth=2pt]{-}(6,-8)(6,10)(24,10)(24,-8)(6,-8)
\psline[linestyle=dashed,linewidth=2pt]{-}(18,4)(18,22)
\psline[linestyle=dashed,linewidth=2pt]{-}(6,-8)(18,4)(36,4)
\psline[linewidth=2pt]{-}(18,22)(36,22)(36,4)(24,-8)
\psline[linewidth=2pt]{-}(18,22)(6,10)
\psline[linewidth=2pt]{-}(36,22)(24,10)
% here is the + bc
%\psline[linestyle=dashed,linewidth=0.25pt]{-}(6,1)(6,-8)(24,-8)(24,1)
%\psline[linestyle=dashed,linewidth=0.25pt]{-}(18,13)(18,4)
%\psline[linestyle=dashed,linewidth=0.25pt]{-}(18,4)(36,4)(36,13)
%\psline[linestyle=dashed,linewidth=0.25pt]{-}(18,4)(6,-8)
%\psline[linestyle=dashed,linewidth=0.25pt]{-}(36,4)(24,-8)
%\uput[-90](27,22){$B'(\underline i)$}
}}
\endpspicture
} \hskip-20pt \hbox{ \psset{unit=0.1cm} \pspicture(-25,-12)(40,40)
\newgray{darkgray}{0.3}
\newgray{lightgray}{0.8}
\rput{0}(-5,0){ \rput{0}(0,0){
%\psline[linewidth=2pt]{-}(36,13)(24,1)(6,1)
%\psline[linestyle=dashed,linewidth=2pt]{-}(6,1)(18,13)(36,13)
\psline[linewidth=1.5pt]{-}(30,11)(22,3)(12,3)
\psline[linestyle=dashed,linewidth=1.5pt]{-}(12,3)(20,11)(30,11)
\psline[linestyle=dashed,linewidth=1.5pt]{-}(18,4)(20,11)(18,22)
\psline[linewidth=1.5pt]{-}(36,4)(30,11)(36,22)
\psline[linewidth=1.5pt]{-}(24,-8)(22,3)(24,10)
\psline[linewidth=1.5pt]{-}(6,-8)(12,3)(6,10)
%\psdots(21,7)
\psline[linewidth=2pt]{-}(6,1)(6,10)(24,10)(24,-2)
\psline[linestyle=dashed,linewidth=2pt]{-}(18,13)(18,22)
%\psline[linewidth=2pt]{-}(18,22)(36,22)(36,10)
%\psline[linewidth=2pt]{-}(18,22)(6,10)
%\psline[linewidth=2pt]{-}(36,22)(24,10)
\psline[linestyle=dashed,linewidth=2pt]{-}(6,1)(6,-8)(24,-8)(24,1)
\psline[linestyle=dashed,linewidth=2pt]{-}(18,13)(18,4)
\psline[linestyle=dashed,linewidth=2pt]{-}(18,4)(36,4)(36,13)
\psline[linestyle=dashed,linewidth=2pt]{-}(18,4)(6,-8)
\psline[linestyle=dashed,linewidth=2pt]{-}(36,4)(24,-8)
\psline[linewidth=2pt]{-}(6,-8)(6,10)(24,10)(24,-8)(6,-8)
%\psline[linestyle=dashed,linewidth=2pt]{-}(18,4)(18,22)
\psline[linestyle=dashed,linewidth=2pt]{-}(6,-8)(18,4)(36,4)
\psline[linewidth=2pt]{-}(18,22)(36,22)(36,4)(24,-8)
\psline[linewidth=2pt]{-}(18,22)(6,10)
\psline[linewidth=2pt]{-}(36,22)(24,10) }}
\endpspicture
} }
%\vskip-10pt
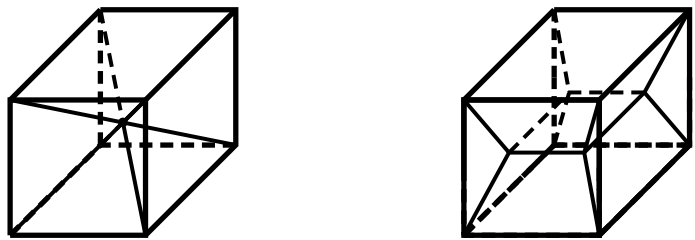
\centerline{figure 1}
%\end{figure}
\smallskip
Consider the Potts model with $q=6$ colors (states)  in a three
dimensional box with boundary condition $i$ on the $i$-th face of
the box. Naively, one might expect that all phases will try to
occupy the region closest to the corresponding piece of the
boundary, which would lead to a phase partition consisting of
symmetric and pyramid-like regions, as can be seen in figure 1,
left. However, at least in the case when the surface tension is
isotropic (which is presumably the case in the limit $T\uparrow
T_c$), there exists a better configuration with lower total
surface free energy. Recall that in this case our desired
interface is simply a minimal surface spanned by the edges of the
box. A picture of the well known solution to this problem can be
seen in figure 1, right. In order to be able to discuss this
example at temperatures $0<T<T_c$, we have to make certain
assumptions about the surface tension $\tau$. We assume that the
sharp simplex inequality holds, that the value of $\tau$ is
minimal in axis directions and that~$\tau$ increases as the normal
vector moves from say $(0,0,1)$ to $(1,1,1)$. (Although these
assumptions are very plausible, none of them has been proved in
dimensions $d\ge 3$). Under these hypotheses, we conjecture that
the phase partition at moderate subcritical temperatures looks
like in figure 2, left. In the limit $T\downarrow 0$, only two
phases survive, as shown in figure 2, right. At $T=0$, there is no
reason for the middle plane to stay centered, in fact, any
horizontal plane is equally likely.
\medskip
%\begin{figure}[bt]
\vskip-40pt \centerline{ \hskip-30pt \hbox{ \psset{unit=0.1cm}
\pspicture(-25,-12)(40,40)
\newgray{darkgray}{0.3}
\newgray{lightgray}{0.8}
\rput{0}(-5,0){ \rput{0}(0,0){
%\psline[linewidth=2pt]{-}(36,13)(24,1)(6,1)
\psline[linewidth=1.5pt]{-}(33,12)(22,1)(9,1)
\psline[linestyle=dashed,linewidth=1.5pt]{-}(9,1)(20,12)(33,12)
\psline[linestyle=dashed,linewidth=1.5pt]{-}(18,4)(20,12)(18,22)
\psline[linewidth=1.5pt]{-}(36,4)(33,12)(36,22)
\psline[linewidth=1.5pt]{-}(24,-8)(22,1)(24,10)
\psline[linewidth=1.5pt]{-}(6,-8)(9,1)(6,10)
%\psdots(21,7)
\psline[linewidth=2pt]{-}(6,1)(6,10)(24,10)(24,-2)
\psline[linestyle=dashed,linewidth=2pt]{-}(18,13)(18,22)
%\psline[linewidth=2pt]{-}(18,22)(36,22)(36,10)
%\psline[linewidth=2pt]{-}(18,22)(6,10)
%\psline[linewidth=2pt]{-}(36,22)(24,10)
\psline[linestyle=dashed,linewidth=2pt]{-}(6,1)(6,-8)(24,-8)(24,1)
\psline[linestyle=dashed,linewidth=2pt]{-}(18,13)(18,4)
\psline[linestyle=dashed,linewidth=2pt]{-}(18,4)(36,4)(36,13)
\psline[linestyle=dashed,linewidth=2pt]{-}(18,4)(6,-8)
\psline[linestyle=dashed,linewidth=2pt]{-}(36,4)(24,-8)
\psline[linewidth=2pt]{-}(6,-8)(6,10)(24,10)(24,-8)(6,-8)
%\psline[linestyle=dashed,linewidth=2pt]{-}(18,4)(18,22)
\psline[linestyle=dashed,linewidth=2pt]{-}(6,-8)(18,4)(36,4)
\psline[linewidth=2pt]{-}(18,22)(36,22)(36,4)(24,-8)
\psline[linewidth=2pt]{-}(18,22)(6,10)
\psline[linewidth=2pt]{-}(36,22)(24,10) }}
\endpspicture
} \hskip-20pt \hbox{ \psset{unit=0.1cm} \pspicture(-25,-12)(40,40)
\newgray{darkgray}{0.3}
\newgray{lightgray}{0.8}
\rput{0}(-5,0){ \rput{0}(0,0){
\psline[linewidth=2pt]{-}(36,13)(24,1)(6,1)
\psline[linestyle=dashed,linewidth=2pt]{-}(6,1)(18,13)(36,13)
%\psdots(21,7)
\psline[linewidth=2pt]{-}(6,1)(6,10)(24,10)(24,-2)
\psline[linestyle=dashed,linewidth=2pt]{-}(18,13)(18,22)
%\psline[linewidth=2pt]{-}(18,22)(36,22)(36,10)
%\psline[linewidth=2pt]{-}(18,22)(6,10)
%\psline[linewidth=2pt]{-}(36,22)(24,10)
\psline[linestyle=dashed,linewidth=2pt]{-}(6,1)(6,-8)(24,-8)(24,1)
\psline[linestyle=dashed,linewidth=2pt]{-}(18,13)(18,4)
\psline[linestyle=dashed,linewidth=2pt]{-}(18,4)(36,4)(36,13)
\psline[linestyle=dashed,linewidth=2pt]{-}(18,4)(6,-8)
\psline[linestyle=dashed,linewidth=2pt]{-}(36,4)(24,-8)
\psline[linewidth=2pt]{-}(6,-8)(6,10)(24,10)(24,-8)(6,-8)
%\psline[linestyle=dashed,linewidth=2pt]{-}(18,4)(18,22)
\psline[linestyle=dashed,linewidth=2pt]{-}(6,-8)(18,4)(36,4)
\psline[linewidth=2pt]{-}(18,22)(36,22)(36,4)(24,-8)
\psline[linewidth=2pt]{-}(18,22)(6,10)
\psline[linewidth=2pt]{-}(36,22)(24,10) }}
\endpspicture
} }
%\vskip-10pt
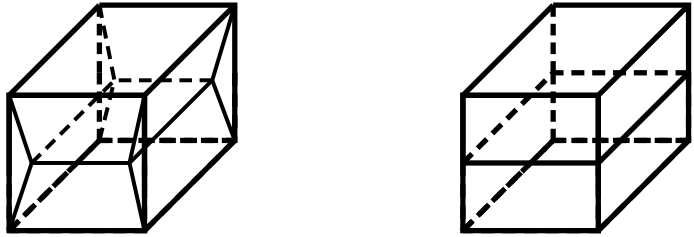
\centerline{figure 2}
%\end{figure}
\smallskip

In the next example we consider the three dimensional Ising model
with free boundary conditions below $T_c$, conditioned on the
event that the average magnetization is positive and does not
exceed $m^*-\eps$, where $\eps$ is a sufficiently small positive
number and $m^*$ denotes the spontaneous magnetization. It is
natural to conjecture that a minimizer of the corresponding
variational problem is a droplet attached symmetrically to one of
the corners of the box. \comment{
%\vskip-15pt
\medskip
%\begin{figure}[ht]
\centerline{ \epsfxsize 8.5cm \epsfysize 6cm \epsfbox{cw.ps} }
\smallskip
\centerline{figure 4}
%\end{figure}
\goodbreak

%\medskip
\begin{figure}[ht]
\centerline{ \epsfxsize 4cm \epsfysize 5cm \epsfbox{dbu.eps} }
\smallskip
\centerline{figure 5}
\end{figure}
} %%%%%%%%%%%%%%%%%%%%%%%%%%%%%%% comment

The single bubble sitting in one of the corners is filled with the
minus phase and in the rest of the box we see the plus phase. The
size of the bubble is determined by $\eps$ and its internal
boundary coincides with the corresponding piece of the surface of
the Wulff crystal.

Another  Wulff-type problem arises by conditioning the $q$-states
Potts model (with say $q\ge 4$) to have a moderate excess of
colors 2 and 3 while imposing 1-boundary conditions on the entire
box. In this case it is conceivable that a so-called ``double
bubble'' is created, consisting of two adjacent macroscopic
droplets filled with the (pure) phases 2 and 3, respectively. The
double bubble is swimming in the phase 1 which fills the rest of
the box. Of course, we might have an excess of color 4 as well; in
this case a further bubble will presumably appear which will be
attached to the previous two bubbles.

For related variational questions concerning soap films and
immiscible fluids, see \cite{LM}.

In fact, by studying questions concerning phase boundaries we are
very quickly confronted with the theory of minimal surfaces, such
as the Plateau problem, corresponding to anisotropic surface
measures. Let $\Om$ be a bounded open set in $\R^3$ with smooth
boundary and let $\ga$ be a Jordan curve drawn on $\bd\Om$ which
separates $\bd\Om$ into two disjoint relatively open sets $\Ga^+$
and $\Ga^-$. Typical configurations in the Ising model on a fine
grid in $\Om$ with plus b.c.s on $\Ga^+$  and minus  b.c.s on
$\Ga^-$ will exhibit two phases separated with an interface close
to a  minimal surface which is a global solution of the following
Plateau type problem:
\[
\mb{minimize } \int_{S} \tau(\nu_S(x))\, d\hdu(x)\ : \mb{ $S$ is a
surface in $\Om$ spanned by $\ga$}
\]
where $\nu_S(x)$ is the normal vector to $S$ at $x$. We remark
that it is conjectured that, as the temperature approaches $T_c$
from below, the surface tension $\tau$ becomes more and more
isotropic and it is conceivable that the solution of the above
minimization problem approaches the solution of the classical
(isotropic) Plateau problem.

\Subsection{Further background}

\vskip-5mm \hspace{5mm}

%{\bf FK percolation and the FK representation of the Potts model.}
There is a beautiful and extremely useful way to decompose the
Ising-Potts model into a certain bond percolation model, called
FK-percolation and some simple 'coloring' procedure discovered by
Fortuin and Kasteleyn \cite{FK}. Consider the $q$-state Potts
model with mixed boundary conditions $\ga$ in a finite  lattice
box $B$ at inverse temperature $\beta$ and set $p = 1-e^{-\beta}$.
Consider a bond percolation model (called FK percolation)
specified by the following formula:
 \[
\fk{\ga,\beta,q}{B}{\eta} = P_p[\eta] \, q^{\#(\eta)} /
Z^{\ga,\beta,q}_B
\]
where $\eta$ is a bond configuration with the property that there
is no open connection between differently colored boundary parts,
$P_p$ is the usual Bernoulli measure with parameter $p$,
$\#(\eta)$ denotes the number of clusters (conn. components) in
the configuration $\eta$ with the rule that identically colored
boundary parts (and their connected components) count as one
single cluster. Finally $Z^{\ga,\beta,q}_B$ is the appropriate
normalizing constant.

In the second step we assign colors to every cluster (and their
sites) as follows: the boundary pieces inherit the color of the
boundary condition, the remaining clusters will get one of the
colors $1,2,...,q$ with probability $1/q$ each independently from
each other. The distribution of the coloring of the sites
corresponds exactly to the Potts model. Note that in the case of
free or constant b.c.s. there is no constraint prohibiting open
connections, and indeed these measures, denoted by
$\fko{f,\beta,q}{B}$ (free b.c.s.) and $\fko{w,\beta,q}{B}$
(``wired'' b.c.s.) behave very similar to regular Bernoulli
percolation. Their thermodynamic limits $\fko{f,\beta,q}{\infty}$
and $\fko{w,\beta,q}{\infty}$ exist as the box size tends to
infinity and we can define the percolation probability as usual
$\theta^*(p) = \fk{*,p}{\infty}{0\conn \infty}$, for $*$ free or
wired. It is easy to check that the order parameter
$\theta(\beta)$
 of the Ising-Potts model agrees with $\theta^w(\beta)$, and
correspondingly the FK model exhibits a percolation phase
transition. Further it is known \cite{G.2}, that $\theta^w(\beta)
= \theta^f(\beta)$ for all but at most countably many values of
$\beta$ and it is conjectured that the equality is valid for all
values except possibly the critical point $\beta_c = \beta_c(q)$.
Moreover this condition is equivalent to the equality of the
thermodynamical limits: $\fko{w,\beta,q}{\infty} =
\fko{f,\beta,q}{\infty}$. We define the set of ``regular'' inverse
temperatures by
\[
\calU(q) = \set{\beta >0}{\theta^w(\beta) =\theta^f(\beta)}.
\]
Although the status of the bonds are dependent, their correlation
tends fast to zero as the distance between them becomes large and
this holds for both the sub and supercritical phase of FK
percolation. This property is of crucial importance in our large
deviation analysis since in the original Potts model there exists
no corresponding asymptotic independence when $\beta > \beta_c$.

Our results are valid above the so called {\em slab-threshold}
$\bchat=\bchat(q,d)$, introduced in \cite{Pi}. This threshold is
conjectured to agree with the critical point and at least in the
case of percolation ($q=1$) this have been proved by Grimmett and
Marstrand \cite{G-M}. It is possible to characterize this
threshold as the smallest value such that when $\beta$ exceeds it,
it is possible to find $\alpha >0$ and $L\ge 1$ such that at least
in the center of the of slabs $S(L,n) = [-L,L]\times
[-n,n]^{d-1}\cap \Zd$ there is ``uniform long range order'', i.e.,
\[
\inf_{n\ge 1}\ \inf_{x,y \in S(L,\alpha n)} \FK{f,
\beta}{S(L,n)}{x\conn y} > 0.
\]
It has been proved in \cite{Pi} that above $\bchat$, $\alpha$
always can be chosen to be one, guaranteeing a strictly positive
probability (uniformly in $n$) for connections within a
sufficiently (depending on $\beta$) thick slab. This property is
crucial for establishing the basic properties of supercritical FK
percolation; the existence of a unique crossing cluster in a box,
its omnipresence, the concentration of its density around the
percolation probability, the exponential tail decay of the
diameter of other clusters in the box, etc. These properties are
then used to establish a renormalization scheme which is essential
for the large deviation analysis of this and the Ising-Potts
model.

\Section{The results}

\vskip-5mm \hspace{5mm}

\Subsection{Historical remarks}

\vskip-5mm \hspace{5mm}

Before we start with the presentation of our results we give a
brief summary of the previous work on this subject. As we have
already mentioned, large deviations theory plays an important role
in this context and not surprisingly the first efforts were
devoted to the study of large deviations of the empirical
magnetization in the Ising model, i.e., the average value of the
spins in a large box. A {\it volume order} large deviation
principle (LDP) has been established for the Ising model by
various authors: Comets, Ellis, F\"ollmer, Orey, Olla \cite{Com,
E, Foe-O, Ol}. The corresponding rate function has been found to
vanish in $[-m^*,m^*]$ where $m^*$ denotes the spontaneous
magnetization. In fact,  it was suspected that the correct order
of decay is exponential to  {\it surface order}. Indeed, Schonmann
\cite{S} found a proof of this conjecture,  valid for any
dimensions and low enough temperatures and Chayes, Chayes and
Schonmann extended the result for the supercritical
$\beta>\beta_c$ regime in the two dimensional case.  F\"ollmer and
Ort \cite{Foe-Ort} investigated  this phenomenon on the level of
empirical measures. Finally, inspired  by the work of Kesten and
Zhang \cite{KZ} on related questions in percolation, Pisztora
\cite{Pi} established surface order upper bounds for the remaining
dimensions $d\ge 3$ above the slab-threshold $\bhc$, introduced in
the same work, which is conjectured to agree with the critical
point~$\beta_c$. In that work a renormalization scheme has been
developed for supercritical Fortuin-Kasteleyn percolation (or
random cluster model) in conjunction with a stochastic domination
argument (generalized and improved in \cite{LSS}) which allows to
control the renormalized process, and so, the original one.

The monograph of Dobrushin, Koteck\'y and Shlosman \cite{DKS}
opened the way to the rigorous study of the phase separation
phenomenon creating thereby an immense interest and activity which
lasts up to the present time. Their analysis, which provided the
first mathematical proof of phase separation, had been performed
in the context of the Ising model. The main tool of their work is
the cluster expansion, which, on the one hand allowed the
derivation of results much finer than necessary to verify the
Wulff construction, on the other hand it restricted the validity
of the results to two dimensions and low temperatures. Significant
improvements of these results in two dimensions have been derived
by Pfister \cite{Pf}, Alexander, Chayes and Chayes \cite{ACC}
(treating percolation), Alexander \cite{Al}, Ioffe \cite{Io1,Io2}.
Finally Ioffe and Schonmann \cite{IOSC} extended the  results of
\cite{DKS} up to $T_c$.

The next challenge was to analyze phase separation for short range
models in higher dimensions. The additional difficulties came
mainly from two sources. First, new techniques have to be
developed to avoid the use of perturbative methods (such as the
cluster expansion) which severely limit the applicability of the
arguments and methods which are specific to two dimensions only
(duality). Second, the emerging geometry is far more complex than
in two dimensions and this requires the use of new tools and
ideas. The complexity of the geometry causes problems also within
the probabilistic analysis (for instance the lack of the skeleton
technique for surfaces) and even the correct formulation of the
results is far from obvious (``hairs'').

The first issue has been resolved by the application of the
aforementioned renormalization technology from \cite{Pi}.
Renormalization arguments lie at the heart of the proof of much of
the intermediate steps (for instance exponential tightness,
decoupling) and even in the remaining parts they play an important
role usually in combination with geometric arguments (interface
lemma, etc.).

To handle the geometric difficulties, the use of appropriate tools
from geometric measure theory has been introduced in the works of
Alberti, Bellettini, Bodineau, Butt\`a, Cassandro, Presutti
\cite{ABCP, BBBP, BBP} and by Cerf \cite{Ce}. In these works also
a novel and very general large deviation framework have been
proposed to tackle the problem. In fact, this framework turned out
to be crucial for the success of the entire approach. It is the
work of Cerf \cite{Ce} in which the first complete analysis of
phase separation in a three dimensional model have been achieved,
namely the asymptotic analysis of the shape of a large finite
cluster (Wulff problem) in percolation.

The results presented in this article have been derived in the
works \cite{Pi, Ce, CP, CP2}. It should be mentioned that in an
independent work \cite{BO} Bodineau carried out an analysis of the
Wulff problem in the Ising model with conclusions slightly weaker
than the results appearing in \cite{CP}.

Finally, for current developments in the field we refer the reader
to the preprint \cite{BP} and the references therein.

\Subsection{Statement of the  results}

\vskip-5mm \hspace{5mm}

{\em Range of validity of the results.} Our results for the
Ising-Potts models hold in the region: $d\ge 3,\
q\in\N\setminus\{0,1\},\ \beta>\bchat(q,d), \ \beta\in
\calU(q,d)$.

At this point it is natural to comment on the case of two
dimensions. Although most of our results should hold for $d=2$,
there are several points in the proofs which would require a
significant change, making the proofs even longer. The main
reason, however,  for not to treat the  two dimensional case  is
that the natural topology for the LDP-s in $d=2$ is not the one we
use (which is based on the distance $\dist$) but a topology based
on the Hausdorff distance.

{\em Surface tension.} From FK percolation  we can  extract a
direction dependent surface tension $\tau(\nu) = \tau(p,q,d,\nu)$,
cf \cite{CP}. For a unit vector~$\nu$, let $A$ be a unit
hypersquare orthogonal to $\nu$, let $\cyl A$ be the cylinder
$A+\R\nu$, then $\tau(\nu)$ is equal to the limit
\[
\lim_{n\to\infty} \! - \LDn\! \Phi_\infty^{p,q}\!\left( \lower 18pt\vbox{ \hbox{inside $n\cyl A$ there exists a
finite set of closed edges $E$} \hbox{cutting $n\cyl A$ in at least 2 unbounded components and} \hbox{the edges of
$E$ at distance less than $2d$ from the boundary} \hbox{of $n\cyl A$ are at distance less than~$2d$
from~$nA$}}\kern-0.5pt \right)\kern-0.7pt
\]
The function~$\tau$ satisfies the weak simplex inequality, is
continuous, uniformly bounded away from zero and infinity and
invariant under the isometries which leave~$\Zd$ invariant (see
section 4 in \cite{CP} for details).

{\em Identification of the phases.} The typical picture which
emerges from the Potts model with mixed b.c.s. at the macroscopic
level is a partition of $\Om$ in maximal $q$ phases corresponding
to the dominant color in that phase. The individual phases need
not be connected. In order to identify the phases we choose first
a sequence of {\em test events} which we regard as characteristic
for that phase. More specificly, for $j=1,2,...,q$ we select
events $E^{(j)}_n$ defined on a $n\times n$ lattice box such that
$\Mu{(j),\beta,q}{\infty}{E^{(j)}_n} \to 1$ as $n\to \infty$. We
may also assume that $E^{(j)}_n\cap E^{(j)}_n = \emptyset$ for
$j\neq i$. For instance, one natural choice is to require that the
densities of the different colors in the box do not deviate more
than some small fraction from their expected value. This will
guarantee that the right mixture of colors occurs which is typical
for that particular pure phase. Alternately, we may request that
the empirical measure defined by the given configuration is close
to the restriction of $\Muo{(j),\beta,q}{\infty}$ to $\La(n)$ with
respect some appropriate distance between probability measures,
etc.

For $x\in\Rd$ and $r>0$ we define the box $\La(x,r)$ by
\[
\La(x,r) = \Set{y\in\Rd}{-r/2 < y_i - x_i \le r/2, \ i=1,\dots,d}
\]
and we introduce an {\it intermediate length scale} represented by
a fixed function $f:\N\to\N$ satisfying \bea\lb{fcond} \limn
n/f(n)^{d-1} = \limn f(n)/\log n = \infty. \eea Given a
configuration in $\Omn$, we say that the point $x\in \Om$ belongs
to the phase $j$, if the event $E^{(j)}_n$ occurs in the box
$\La(x, f(n)/n)$. For $j=1,2,..,q$, we denote by $A^{(j)}_n$ the
set of points in $\Om$  belonging to the phase $j$ and set
$A^{(0)}_n = \Om\setminus \cup_{j}A^{(j)}_n$ (indefinite phase).

The random partition of $\Om$, $\vec{A}_n =
(A^0_n,A^1_n,\dots,A^q_n)$, is called the {\it empirical phase
partition}. Our first result shows that up to super-surface order
large deviations, the region of indefinite phase $\An{0}$ has
negligible density, i.e., the pure phases fill out the entire
volume.
\begin{theorem}\lb{thmspin}
Let $d\ge 3$, $q\in\N\setminus\{0,1\}$, $\beta >\bchat$, $\beta\in
\calU(q,d)$. For $\de >0$, \bean \lsn \LDn \Mu{}{n}{\vol(\An{0})
>\de }= - \infty. \eean
\end{theorem}

Although Theorem \ref{thmspin} guarantees that the pure phases
fill out the entire volume (up to negligible density) but it does
not exclude the possibility that the connected components of the
pure phases are very small. For instance, they could have a
diameter not much larger than our fixed intermediate scale (in
which case they would be invisible on a macroscopic scale.) If
this happened, the total area of the phase boundaries would be
exceedingly high. Before stating our next result, which will
exclude this possibility, we introduce some geometric tools.

We define a (pseudo) metric, denoted by $\dist$, on the set
$\calB(\Om)$ of the Borel subsets of $\Om$ by setting
\bea\lb{distdef} \forall A_1,A_2\in\calB(\Om)\qquad
\dist(A_1,A_2)\,=\,\vol(A_1\Delta A_2). \eea We consider then the
{\em space of phase partitions} $P(\Om,q)$ consisting of
$q+1$-tuples $(A^0, A^1, \dots, A^q)$ of Borel subsets of $\Om$
forming a partition of $\Om$. We endow $P(\Om,q)$ with the
following metric:
\[
\DistP\Big( (A^0, \dots, A^q), (B^0, \dots, B^q)\Big) =
\sum_{i=0,\dots,q}\dist(A^i,B^i).
\]
In order to define the  surface energy  $\calI$ of a phase
partition $\vec{A_n}$, we recall some notions and facts from the
theory of sets of finite perimeter, introduced initially by
Caccioppoli and subsequently developed by De Giorgi, see for
instance \cite{GI,EVGA}. The perimeter of a Borel set~$E$ of
$\R^d$ is defined as
$$\calP(E)\,=\,
\sup\,\Big\{\, \int_E\text{div}\,f(x)\,dx: f\in
C^\infty_0(\R^d,B(1))\,\Big\}\,$$ where $C^\infty_0(\R^d,B(1))$ is
the set of the compactly supported $C^\infty$ vector functions
from $\R^d$ to the unit ball $B(1)$ and $\text{div}$ is the usual
divergence operator. The set~$E$ is of finite perimeter if
$\calP(E)$ is finite. A unit vector~$\nu$ is called the measure
theoretic exterior normal to~$E$ at~$x$ if
$$\lim_{r\to 0} r^{-d}\vol(B_-(x,r,\nu)\setminus E)\,=\,0\,,\quad
\lim_{r\to 0} r^{-d}\vol(B_+(x,r,\nu)\cap E)\,=\,0\,.$$ Let $E$ be
a set of finite perimeter. Then there exists a certain subset of
the topological boundary of $E$, called the {\em reduced
boundary}, denoted by $\rb E$, with the same $d-1$ dimensional
Hausdorff measure as $\partial E$, such that at each $x\in \rb E$
there is a measure theoretic exterior normal to~$E$ at~$x$. For
practical (measure theoretic) purposes, the reduced boundary
represents the boundary of any set of finite perimeter, for
instance,
 the following generalization of Gauss Theorem holds:
For any vector function~$f$ in~$C^1_0(\R^d,\R^d)$,
$$\int_E\text{div}\, f(x)\, dx\,=\,
\int_{\rb E}f(x)\cdot\nu_E(x)\, \h^{d-1}(dx)\,.$$ (For more on
this see e.g. the appendix in \cite{CP} and the references there.)

The {\em surface energy} $\calI$ of a phase partition $(\An{0}
\An{1}, \dots, \An{q})\in P(\Om,q)$ is defined as follows:
\begin{itemize}
\item[-]
for any $(A^0, A^1, \dots, A^q)$ such that either $A^0\nempty$ or
one set among $A^1, \dots, A^q$ has not finite perimeter, we set
$\calI(A^0, \dots, A^q) = \infty,$
\item[-]
for any $(A^0, A^1, \dots, A^q)$ with $A^0\empty$ and  $A^1,
\dots, A^q$ having finite perimeter, we set
\end{itemize}
\begin{eqnarray*}
  \calI(A^0, \dots, A^q) & = & \sum_{i=1,\dots,q} \rec{2}\int_{\bd^*A^i\cap\Om}\tau(\nu_{A_i}(x))\,
  d\calH^{d-1}(x)\\
  & & + \sum_{i, j =1,\dots,q \atop i\neq j} \int_{\bd^*A^i\cap \Ga^j}
  \tau(\nu_{A_i}(x))\,  d\calH^{d-1}(x).
\end{eqnarray*}
%\bean
%\hspace{-1cm} \calI(A^0, \dots, A^q)& =& \sum_{i=1,\dots,q}
%\rec{2}\int_{\bd^*A^i\cap\Om}\tau(\nu_{A_i}(x))\,  d\calH^{d-1}(x) \\
%&& \hspace{4cm}  +
%\sum_{i, j =1,\dots,q \atop i\neq j}
%\int_{\bd^*A^i\cap \Ga^j}
%\tau(\nu_{A_i}(x))\,  d\calH^{d-1}(x)
%\eean
Note that $\calI$ depends on $\tau$ and the boundary conditions
$\ga = (\Ga^1,...,\Ga^q)$. The first term in the above formula
corresponds to the interfaces present in $\Om$, while the second
term corresponds to the interfaces between the elements of the
phase partition and the boundary $\Ga$. It is natural to define
the {\em perimeter of the phase partition} by using the same
formula with $\tau$ replaced by constant one. Since $\tau$ is
uniformly bounded away from zero and infinity at any temperature,
the surface energy can be bounded by a multiple of the perimeter
and vice versa. It is known that the surface energy $\calI$ and
the perimeter~$\calP$ are lower semi continuous  and their level
sets of the form $\calI^{-1}[0, K]$ are compact  on the space
$(\calB(\Rd),\dist)$.

The next result states that up to surface order large deviations
(and the constant can be made arbitrarily large by adjusting the
bound $K$ below) the empirical phase partition will be
(arbitrarily) close to the set of phase partitions with perimeter
not exceeding $K$.
\begin{theorem}\lb{exptight}
Let $d\ge 3$, $q\in\N\setminus\{0,1\}$, $\beta >\bchat$, $\beta\in
\calU(q,d)$. For $\de >0$,
\[
\lsn \LDn \MU{}{n}{ \DistP(\vec{A}_n, \, \calI^{-1}[0,K])> \delta}
\le - c\, K.
\]
\end{theorem}

Our fundamental result is a large deviation principle (LDP) for
the empirical phase partition $(\An{0}, \An{1}, \dots, \An{q})$.
\begin{theorem}\lb{spinLDP}
The sequence $\big(\vec{A}_n\big)_{n\in\N} =
\big((\An{0},\An{1},\dots,\An{q})\big)_{n\in\N}$ of the empirical
phase partitions of $\Om$ satisfies a LDP in $(P(\Om,q),\DistP)$
with respect to $\mu_n$ with speed $n^{d-1}$ and rate function
$\calI-\min_{P(\Om,q)}\calI$, i.e., for any Borel subset $\Bbb E$
of $P(\Om,q)$,
\begin{eqnarray*}
  - \inf_{\cdero} \,  \calI +\min_{P(\Om,q)}\calI & \le & \lin \LDn \mu_n\Big[{\vec{A}_n\in  {\Bbb E}}\Big] \\
  & \le & \lsn \LDn \mu_n\Big[{ \vec{A}_n\in {\Bbb E}}\Big] \\
  & \le & - \inf_{\bar{\Bbb E}} \,  \calI +\min_{P(\Om,q)}\calI.
\end{eqnarray*}
\end{theorem}

Note that the constant $\min_{P(\Om,q)}\calI$ is always finite.
Every large deviation result includes a (weak) law of large
numbers; here the corresponding statement is as follows:
Asymptotically, the empirical phase partition will be concentrated
in an (arbitrarily) small neighborhood of the set of partitions
minimizing the surface energy. In other words, on the macroscopic
level, the typical phase partition will coincide with some of the
minimizers of the variational problem, in agreement with the
phenomenological prediction. Of course, the LDP states much more
than this, in particular we will be able to extract similar
statements for restricted ensembles. Recall that imposing mixed
boundary conditions is not the only way to force the system  to
exhibit coexisting phases. In the Wulff problem in the Ising model
context, for instance, a restricted ensemble is studied which is
characterized by an artificial excess of say minus spins in the
plus phase. Technically this can be achieved by conditioning the
system to have a magnetization larger than the spontaneous
magnetization while imposing plus b.c.s.

The  next result describes the large deviation behavior of the
phase partition in a large class of restricted ensembles. Although
it is a rather straightforward generalization of Theorem
\ref{spinLDP}, we state it separately because of its physical
relevance.

Let $(G_n)_{n\ge 1}$ be a sequence of events, i.e., sets of spin
configurations, satisfying the following two conditions: first
there exists a Borel subset $\GGG$ of $P(\Om,q)$ such that the
sequence of events $(G_n)_{n\in\N}$ and $(\{\vec{A}_n\in
\GGG\})_{n\in\N}$ are exponentially equivalent, i.e., \bea
%&&\ \mb{ there exists a Borel subset $\GGG$ of $P(\Om,q)$ such that the
%sequence of events
%$(G_n)_{n\in\N}$ and $(\{\vec{A}_n\in \GGG\})_{n\in\N}$
%are exponentially equivalent, i.e.}
%\nonumber \\
\lsn \LDn \mu_n\Big[ \symdif{G_n}{ \{\vec{A}_n\in \GGG \}}\Big] =
-\infty, \lb{G1} \eea where $\triangle $ denotes the symmetric
difference. Second, the following limit exists and is finite:
\bea\lb{G2} \calI_G = \limn \LDn \mu_n[\, G_n ]
>-\infty.
\eea The sequence of events $(G_n)_{n\ge 1}$ determines a
restricted (conditional) ensemble. \comment{ For instance, in the
case of the Wulff problem in the Ising model,
$$G_n = \Big\{ |\Omn|^{-1}\sum_{x\in \Omn} \si(x) < m\Big\}$$
where $m$ is a fixed value strictly smaller than the spontaneous
magnetization, but sufficiently close to it to ensure that the
adequately dilated Wulff crystal fits into the region $\Om$.}
%satisfying $m<m^*$ and $m^*$ denotes the
%$m$ is a fixed value satisfying $m<m^*$ and $m^*$ denotes the
%spontaneous magnetization.
Note that if \bea\lb{G3} \inf\limits_{\bgcirc} \,  \calI =
\inf\limits_{\bar{\Bbb G}} \,  \calI > -\infty, \eea then
%the large deviation principle of
Theorem \ref{spinLDP} implies that \bref{G2} is satisfied, with
$\calI_G= \inf\limits_{{\Bbb G}} \,  \calI$.
\begin{theorem}\lb{condspinLDP}
Assume that the sequence $(G_n)_{n\ge 1}$ satisfies \bref{G1} and
\bref{G2} and define for each $n\ge 1$ the conditional measures
\[
\mu^G_n = \mu_n(\, \cdotp | G_n).
\]
Then the sequence $\big(\vec{A}_n\big)_{n\ge 1}$ of the empirical phase partitions of $\Om$ satisfies a LDP in
$(P(\Om,q),\DistP)$ with respect to $\mu^G_n$ with speed $n^{d-1}$ and rate function $\calI- \calI_G$, i.e., for
any Borel subset $\Bbb E$ of $P(\Om,q)$,
\begin{eqnarray*}
  - \inf_{\begcirc} \,  \calI + \calI_G & \le & \lin \LDn \mu^G_n\Big[{\vec{A}_n\in  {\Bbb E}}\Big] \\
  & \le & \lsn \LDn \mu^G_n\Big[{ \vec{A}_n\in {\Bbb E}}\Big] \\
  & \le & - \inf_{\overline{{\Bbb E}\cap\GGG }} \,  \calI +\calI_G.
\end{eqnarray*}
\end{theorem}

Theorem~\ref{condspinLDP} gives a rigorous verification of the
basic assumption underlying the phenomenological theory, namely,
that in a given ensemble, the typical configurations are those
minimizing the surface free energy. A general compactness argument
implies the existence of at least one such minimizer. However, in
most examples one cannot say much about the minimizers themselves.
(One notable exception is the Wulff problem.) The difficulty stems
from the fact that the surface tension~$\tau$ is anisotropic and
almost no quantitative information about its magnitude is
available. Moreover,  the corresponding variational problems are
extremely hard even in the isotropic case and the (few) resolved
questions represent the state of the art in the calculus of
variations. For instance, a famous conjecture related to the
symmetric double-bubble in the three dimensional case with
isotropic surface energy (perimeter) has only been resolved
recently \cite{HHS} and the asymmetric case remains unresolved
(even in the isotropic case).

%{\bf Application to the Wulff and multiple bubble problem.}
We show next how Theorem~\ref{condspinLDP} can be applied to the
Wulff and multiple bubble problem.
%Let  $q\in\N\setminus\{0,1\}$.
We take pure boundary conditions with color~$1$, that is,
$\Ga^1=\Ga$, $\Ga^2=\dots=\Ga^q=\emptyset$. Let $s_2,\cdots,s_q$
be $q-1$ real numbers larger than or equal to $(1-\theta)/q$. We
set
$$\forall i\in\{2,\dots,q\}\qquad
v_i\,=\,\vol(\Om)\theta^{-1}(s_i-(1-\theta)/q).$$ We define next
the events
$$\forall n\in\N\qquad
G_n\,=\,\{\forall i\in\{2,\dots,q\}\quad\calS_n(i)\geq s_i\}$$ and
the collection of phase partitions
$${\Bbb G}(v_2,\dots,v_q)=
\{\, \vec{A}=(A_0,A_1,\dots,A_q)\in P(\Om,q): \vol(A_2)\geq
v_2,\dots, \vol(A_q)\geq v_q\,\}.$$ It can be shown  that the
sequences of events
$$(G_n)_{n\in\N}
\qquad\text{ and }\qquad (\vec{A}_n\in{\Bbb
G}(v_2,\dots,v_q))_{n\in\N}$$ are exponentially equivalent, i.e.,
they satisfy the condition~\bref{G1}. In order to ensure
condition~\bref{G3}, we suppose that the minimum of the surface
energy $\calI$ over ${\Bbb G}(v_2,\dots,v_q)$ is reached with a
phase partition having no interfaces on the boundary $\Ga$. More
precisely, we suppose that the following assumption is fulfilled.

\noindent{\bf Assumption.} The region $\Om$ and the real numbers
$v_2,\dots,v_q$ are such that there exists
$\vec{A}^*=(A_0^*,A_1^*,\dots,A_q^*)$ in ${\Bbb G}(v_2,\dots,v_q)$
such that
$$\displaylines{
\calI(\vec{A}^*)\,=\, \min\,\{\,\calI(\vec{A});\vec{A}\in {\Bbb
G}(v_2,\dots,v_q)\,\}\cr \forall i\in\{\,2,\dots,q\,\}\qquad
d_2(A_i^*,\Ga)>0. }$$ We expect that this assumption is fulfilled
provided the real numbers $v_2,\dots,v_q$ are sufficiently small
(or equivalently, $s_2,\dots,s_q$ are sufficiently close to
$(1-\theta)/q$), depending on the region $\Om$. This is for
instance the case when $q=2$. Indeed, let $\calW_\tau$ be the
Wulff crystal associated to $\tau$. We know that $\calW_\tau$ is,
up to dilatations and translations, the unique solution to the
anisotropic isoperimetric problem associated to~$\tau$. For $v_2$
sufficiently small, a dilated Wulff crystal
$x_0+\alpha_0\calW_\tau$ of volume $v_2$ fits into $\Om$ without
touching $\Ga$, and the phase partition $\vec{A}^*=(\emptyset,
\Om\setminus (x_0+\alpha_0\calW_\tau), x_0+\alpha_0\calW_\tau)$
answers the problem. In the case $q>2$, we expect that a
minimizing phase partition corresponds to a multiple bubble having
$q-1$ components.

%\begin{lemma}\lb{multbub}
%\end{lemma}

Under the above assumption, we claim that the collection of phase
partitions ${\Bbb G}(v_2,\dots,v_q)$ satisfies~\bref{G3}. For
$\lambda>1$, we define
$$\vec{A}^*(\lambda)=\Big(\emptyset,
\Om\setminus\bigcup_{2\le i\le q}\lambda A_i^*, \lambda
A_2^*,\dots,\lambda A_q^*\Big).$$ Since by hypothesis the sets
$A_2^*,\dots,A_q^*$ are at positive distance from $\Ga$, for
$\lambda$ larger than~$1$ and sufficiently close to~$1$, the phase
partition $\smash{\vec{A}^*(\lambda)}$ satisfies
$$\vec{A}^*(\lambda)\,\in\,{\Bbb G}(\lambda^dv_2,\dots,
\lambda^dv_q)\,\subset\, {\bgcirc}(v_2,\dots, v_q)$$ and moreover
$\calI(\vec{A}^*(\lambda))=\lambda^{d-1}\calI(\vec{A}^*)$. Sending
$\lambda$ to~$1$, and remarking that \linebreak ${\Bbb
G}(v_2,\dots,v_q)$ is closed, we see that ${\Bbb
G}(v_2,\dots,v_q)$ satisfies~\bref{G3}. Thus we can apply
Theorem~\ref{condspinLDP} with the sequence of events
$(G_n)_{n\in\N}$, thereby obtaining a LDP and a weak law of large
numbers for the conditional measures $\mu^G_n = \mu_n(\, \cdotp |
G_n)$. In the particular case $q=2$, we obtain again the main
result (Wulff problem) of our previous paper \cite{CP}. In the
more challenging situations $q> 2$, the unresolved questions
concerning the macroscopic behavior of such systems belong to the
realm of the calculus of variations.

%%%%%%%%%%%%%%%%%%%%%%%%%%%%%%%%%%%%%

\label{lastpage}

\end{document}